\documentclass{amsart}

\newcommand {\uinf} {\|u\|_\infty}
\newcommand {\lam}  {\lambda}
\newcommand {\eps}   {\epsilon}
\newcommand {\real} {\mathbb{R}}
\newcommand {\R}{\real}
\newcommand {\Z}{\mathbb{Z}}
\newcommand {\C}{\mathbb{C}}
\newcommand {\mint} {\dot{M}}
\newcommand {\ra} {\rightarrow}
\newcommand {\rhohat} {\hat{\rho}}
\newcommand {\supp}{\text{supp}\,}
\newcommand {\dist}{\text{dist}\,}
\newcommand {\mtilde}{\tilde{M}}
\newcommand {\kdir}{K^{\text{dir}}}
\newcommand {\distr}{\mathcal{D}'}
\newcommand {\lamone}{\lam^{-1}}
\newcommand {\spcheck}{\check{\makebox[2mm]{ }}}
\newcommand {\ybar}{{\tilde{y}}}
\newcommand {\vol}{{\rm vol}}

\newtheorem{theorem}{Theorem}
\newtheorem{proposition}[theorem]{Proposition}
\newtheorem{lemma}[theorem]{Lemma}
\newtheorem{cor}[theorem]{Corollary}

\begin{document}
\title{Uniform bounds for eigenfunctions of the Laplacian on manifolds 
with boundary}
\author{D.\ Grieser}
\thanks{This work was partially supported by NSF grant  DMS-9306389 and
by the Mathematical Sciences Research Institute, Berkeley}
\address{Humboldt-Universit\"at Berlin, Institut f\"ur Mathematik,
Sitz: Rudower Chaussee 25, 10099 Berlin, Germany\\
e-mail: grieser@mathematik.hu-berlin.de
}

\begin{abstract}
Let $u$ be an eigenfunction of the Laplacian on a compact 
manifold with boundary, with
Dirichlet or Neumann boundary conditions, and let $-\lambda^2$ be the
corresponding 
eigenvalue. 
We consider the problem of estimating $\max_M u$ in terms of $\lambda$,
for large $\lambda$, assuming $\int_M u^2=1$.
We prove that $\max_M u\leq C_M \lambda^{(n-1)/2}$, which is optimal for some
$M$. Our proof simplifies some of the arguments used before for such problems.
We review the 'wave equation method' and discuss some special cases which
may be handled by more direct methods.
\end{abstract}

\maketitle

\section {Introduction}

Let $(M,g)$ be a smooth, compact Riemannian manifold of dimension $n\geq2$,
with smooth boundary $\partial M$.
Let $\Delta$ denote the Laplace-Beltrami operator on functions on $M$.
Consider a solution of the eigenvalue problem, with Dirichlet or Neumann boundary
conditions,
\begin{align} \label{evproblem} 
\begin{split}
 (\Delta + \lam^2)u & = 0,\\
    u_{|\partial M} & = 0 \quad\text{ or }\quad \partial_n u_{|\partial M}=0,
\end{split}
\end{align}
($\partial_n$ denotes the normal derivate)
normalized by the condition
\begin{equation} \label{lzwo}
 \|u\|_2 = 1.
\end{equation}
The subscript $p$ will always indicate the $L^p(M)$ norm.
In this paper we consider the problem of bounding 
$  \|u\|_\infty = \max_{x\in M} |u(x)|$ in terms of $\lambda$,
for large $\lam$. 
The size of $\|u\|_\infty$ may be considered as a rough measure for
how unevenly the function $u$ is distributed over $M$: If $u$ is 'small'
(say $o(1)$ as $\lambda\to\infty$) on a 'large' set $S$ (say, of area 
$|S|=|M|-\delta$)
then the condition \eqref{lzwo} forces $\|u\|_\infty$ to be at least of order
$\delta^{-1/2}$; therefore upper bounds on $\|u\|_\infty$ imply
lower bounds on $\delta$, i.e. the area of the set where $u$ is 'concentrated'.
It may happen that $\delta\to0$ as $\lambda\to\infty$, for certain (sequences
of) eigenfunctions, for example on
the sphere or the disk.
Also, upper bounds on $\|u\|_\infty$ yield upper bounds for the
multiplicities of the eigenvalues of $\Delta$.

In addition to single eigenfunctions we also consider sums of eigenfunctions of
the form
$$ u_I(x) = \sqrt{\sum_{j:\lambda_j\in I} |u_j(x)|^2}$$
for finite  intervals $I\subset \R$. Here $u_1, u_2,\ldots$ is any orthonormal basis
of eigenfunctions with eigenvalues $\lambda_1<\lambda_2\leq\ldots\to\infty$. 
(The sum is independent of the
choice of basis.)
A simple argument using a suitable version of
Sobolev's embedding theorem (see, for example, \cite{Hor:ALPDOIII}, Thm.\ 17.5.3)
shows that 
\begin{equation} \label{sobestimate}
\sup_x u_{[0,\lambda]}(x) \leq C\lam^{n/2}.
\end{equation}
(Here and everywhere $C$ will denote some constant only depending on $(M,g)$.)
Together with  Weyl's law 
\begin{equation}\label{weyl}
\#\{j:\,\lambda_j\leq\lambda\}\sim \gamma_M \lambda^n 
\quad\text{as } \lambda\to\infty
\end{equation}
(with $\gamma_M=(2\pi)^{-n}\vol_{{\rm eucl}}(B^n)\vol(M)$, $B^n$=unit ball in
 $\R^n$) 
 this shows that for each $x$
the average size of $u_1(x),u_2(x),\ldots$ is of order $O(1)$.
However, if $M$ is a sphere or a ball, for example,
 then there is a subsequence $u_{j'}$ of
eigenfunctions with $\|u_{j'}\|_\infty\geq c\lam_{j'}^{(n-1)/2}$, for some constant
$c>0$ (see Section \ref{secdisk}). Our main result shows that this is the worst
possible rate of growth:

\begin{theorem} \label{mainthm}
Let $M^n$ be a compact Riemannian manifold with boundary. There is a constant
$C=C(M)$ such that any solution of \eqref{evproblem} satisfies
\begin{equation} \label{thmbound}
 \uinf \leq C\lam^{(n-1)/2}.
\end{equation}
\end{theorem}
From the theorem one easily derives the following corollary. This is well-known
(see 'Related Results' below).

\begin{cor} \label{cor}
Under the same assumptions as in the theorem, the multiplicity of the
eigenvalue $\lambda^2$ of $-\Delta$ is at most $C \lambda^{n-1}$.
\end{cor}

The theorem is already non-trivial in the case of a domain in $\R^n$
with the Euclidean metric: 'Interior' estimates are simpler in this
case due to the constant coefficients of $\Delta$ (see Section \ref{seceucl}),
but the boundary must still
 be curved, and this causes the main difficulty.
While the bound \eqref{thmbound} is optimal for balls and spheres it
 is far from optimal for a torus
and a rectangle (see Section \ref{sectorus}), 
and possibly for more general $M$ under assumptions
on the curvature. For example, in the case of negative curvature,
the results of
B\'{e}rard \cite{Ber:WECRMWCP} imply a bound of
$C\lambda^{(n-1)/2}/\log\lambda$.
 The first time that the exponent in 
\eqref{thmbound} was improved in any case of nonzero curvature was in the
paper (\cite{IwaSar:LNEAS}) by Iwaniec and Sarnak, for
certain arithmetic surfaces $M$. In the opposite direction, one may ask for
which $M$ the $O(\lambda^{(n-1)/2})$ estimate may not be replaced by
$o(\lambda^{(n-1)/2})$. See \cite{SogZel:RMWMEG} for results in the case without
boundary. 
We do not analyze the dependence of the constant $C$ in 
\eqref{thmbound} on the metric. The methods used imply that it is bounded in
terms of a finite number (depending on $n$) of derivatives of the metric
and of the geodesic curvature of $ \partial M$.
\medskip

{\em Outline of the proof of Theorem \ref{mainthm}.}
 The idea is to use the standard wave kernel method
outside a boundary layer of
width $\lam^{-1}$ and a maximum principle argument inside that layer.

Let us first recall the wave kernel method (cf.\ \cite{Hor:SFEO,Sog:FICA}).
Certain weighted
 sums over many eigenfunctions turn out to be easier to estimate
than single eigenfunctions, since they have a  local
 character. More precisely,
given $\eps>0$,
choose a Schwartz function $\rho$ such that
$$ \rho\geq 0, \quad \rho_{|[0,1]}\geq 1, \quad 
\supp (\hat{\rho}) \subset (-\eps,\eps).$$
Here $\hat{\rho}(t) = \int_{-\infty}^\infty e^{-it\lambda}\rho(\lambda)
\, d\lambda$ denotes the Fourier transform. Proving the 
existence of such a $\rho$
is an easy exercise.

If $u_1,u_2,\ldots$ is an orthonormal basis of real eigenfunctions of the
Dirichlet Laplacian, with eigenvalues $\lam_1<\lam_2\leq\ldots \ra \infty$,
then we consider the sum, convergent by \eqref{sobestimate},
\begin{equation}\label{sum}
 \sum_j \rho(\lam-\lam_j) u_j(x)u_j(y).
\end{equation}
Actually, this is the integral kernel of the operator
$\rho(\lam-\sqrt{-\Delta})$, but we won't use this fact.
If we write, using Fourier's inversion formula,
\begin{align*}
\rho(\lam-\lam_j) & = \frac1{2\pi} \int_{-\infty}^\infty
                                    \rhohat(t) e^{it(\lam-\lam_j)}\,dt\\
                  & = [\rhohat(t)e^{-it\lam_j}]\spcheck(\lam)\\
		  & = 2[\rhohat(t) \cos(t\lam_j)]\spcheck(\lam) 
			  - \rho(\lam+\lam_j), 
\end{align*}
(the superscript $\spcheck$ denotes inverse Fourier transform $t\ra \lam$)
then we get the important identity
\begin{equation} \label{wavemethod}
\sum_j \rho(\lam-\lam_j)u_j(x)u_j(y) = 2[\rhohat(t)K(t,x,y)]\spcheck
			    (\lam) + O(\lam^{-\infty})
\end{equation}
where the error term is small by (\ref{sobestimate}) and rapid decay
of $\rho$, and 
\begin{equation} \label{wavek}
K(t,x,y) = K_M(t,x,y) = \sum_j \cos(t\lam_j) u_j(x)u_j(y)
\end{equation}
is the {\em wave kernel}, 
i.e.\ for each fixed $y\in \mint$ (the interior of $M$)
 it is the solution of
\begin{align}
(\frac{\partial^2}{\partial t^2} - \Delta_x)K(t,x,y) & = 0 \quad 
			       \text{ in } \real_t\times \mint_x \notag\\
 K(0,x,y) & = \delta_y(x)\notag\\
 \frac{\partial}{\partial t}K(t,x,y) &= 0 \label{waveeqn}\\
 K_{|x\in \partial M} & = 0.\notag
\end{align}
The convergence of \eqref{wavek} is in the sense of distributions
(i.e.\ weakly) in $t$, for each fixed $x,y$, as 
the argument leading to \eqref{wavemethod} shows, for example.

The point of \eqref{wavemethod}
 is that $K$ may be analyzed directly from \eqref{waveeqn}.
In particular, solutions of the wave equation have {\em finite propagation
speed} (see \cite{Hor:ALPDOIII}, Lemma 17.5.12, for 
example, for the easy proof using energy
estimates); for $K$ this means
$$ \supp K \subset \{(t,x,y):\, \dist(x,y) \leq |t| \}$$
and that
$K(t,x,y)$ depends only on the data  (i.e.\ $M$ and the metric $g$)
 in $B_t(x,y) := \{z:\, \dist(x,z)
+ \dist(y,z) \leq |t|\}$. Therefore, \eqref{wavemethod} shows that 
the sum \eqref{sum} depends
(up to an error $O(\lambda^{-\infty})$) only on the data in $B_\eps(x,y)$,
and this is the {\em local character} mentioned before.

We will use (\ref{wavemethod}) only on the diagonal, i.e.\ 
for $x=y$. From the assumptions on $\rho$
we get
\begin{equation} \label{uwave}
 u_{[\lam-1,\lam]}(x)^2
 \leq 2[\rhohat(t)K(t,x,x)]\spcheck(\lam) + O(\lam^{-\infty}),
\end{equation}
so the theorem would follow from an estimate
\begin{equation} \label{mainest}
| \int e^{it\lam} \rhohat(t) K(t,x,x)\,dt | \leq C\lam^{n-1}. 
\end{equation}
Since $\rhohat$ is smooth, this is a statement about 
the singularities of $K$ as a distribution in $t$, uniformly
in the parameter $x$.

As a first step one may now obtain the {\em interior estimate}
\begin{equation}\label{int}
u_{[\lam-1,\lam]}(x) \leq C_\eps \lam^{(n-1)/2}\quad\text{ if }
 \dist(x,\partial M)
> \eps.
\end{equation}
In the case of a domain $M\subset \real^n$ one has $K_M=K_{\R^n}$
for $|t|<2\,\dist(x,\partial M)$ by locality (finite propagation speed),
so  \eqref{int} follows
easily from \eqref{uwave} by using
 the explicit expression for the Euclidean wave kernel in terms of 
an $x$-space Fourier transform, see Section \ref{seceucl}. 
In the case of a nonflat metric, $K$ can be analyzed near
$t=0$ by the geometric optics approximation (see \cite{Hor:SFEO},\cite{Hor:ALPDOIII}) 
or the scaling
technique introduced by Melrose in \cite{Mel:TWG}, and this gives
 the interior estimate
in that case.
Another method to obtain \eqref{int}
 uses the Hadamard parametrix for the resolvent, see
\cite{Ava:UEGRM}, \cite{Sog:CLNSCSOEOCM}, for example.

Since $|u(x)|\leq u_{[\lam-1,\lam]}(x)$ for solutions $u$
 of \eqref{evproblem},
\eqref{int} gives in particular
\begin{equation}\tag{\ref{int}'}
|u(x)| \leq C_\eps \lam^{(n-1)/2}\quad\text{ if } x\in M_\eps.
\end{equation}

The claim of Theorem \ref{mainthm} is that $C_\eps$ in (\ref{int}') may be
chosen independent of $\eps$.
When $x$ approaches the boundary in the argument above, 
one has to either shrink
the support of $\rhohat$ -- which means increasing its maximum since
$\int\rhohat = \rho(0) = 1$   (and this gives only $C_\eps \leq 
C\eps^{-1/2}$
in \eqref{int}) --
or analyze the wave kernel in the presence of boundary conditions. 
The difficulty with the latter   is that the geometric optics construction
becomes much more complicated near the boundary because of diffraction
and multiple reflection of geodesics, and a satisfactory parametrix for
the wave equation has so far only been constructed near points where the
boundary is strictly convex or concave (see \cite{MelTay:BPWEGGR}, \cite{Mel:WCMCB}). 
This parametrix
was used in \cite{Gri:LBESPLNCB}\ to prove Theorem \ref{mainthm} near
concave boundary points (e.g.\ near the inner circle of an annular domain).
On the diagonal $x=y$ the singularities
of $K$ have been analyzed in sufficient detail near an arbitrary smooth
boundary by Ivrii, Melrose and
H\"ormander (\cite{Ivr:STSALBOMWB,Mel:TWG,Hor:ALPDOIII}) to allow us to obtain
the desired estimates
 at points $x$ outside a boundary layer of width
$\lamone$.
 In the boundary layer,
a simple maximum principle argument then completes  the proof, because there 
$u$ can have at most one oscillation in the direction perpendicular
to the boundary.
\medskip

{\em Related results.}
Estimates of $\|u\|_\infty$ for large $\lam$ are closely related to asymptotic improvements
over the bound \eqref{sobestimate} on
$u_{[0,\lam]}$ (which is often referred to as the
'spectral function of $\Delta$ on the diagonal'). 
That improvements might be possible is 
suggested by the observation that $\int_M u_I^2 = \# \{j:\,\lam_j \in I\}$,
and Weyl's law \eqref{weyl}.  
Carleman \cite{Car:PAFFMV} was the first to prove the interior pointwise
asymptotics corresponding to \eqref{weyl} (as $\lam\to\infty$)
\begin{equation}\label{carl}
u_{[0,\lam]}(x) = \gamma' \lam^{n/2} + o_\eps (\lam^{n/2}),
\quad\text{for } \dist(x,\partial M)
> \eps
\end{equation}
(for domains in $\R^n$; for manifolds see \cite{MinPle:SPELORM}, \cite{Gar:DPLEPDE}).
Here $\gamma'=(2\pi)^{-n/2}\sqrt{\vol_{\rm{eucl}}(B^n)}$.
The error term was improved to
\begin{equation}\label{avak}
u_{[0,\lam]}(x) = \gamma' \lam^{n/2} + O_\eps(\lam^{(n-1)/2}), 
\quad\text{for }  \dist(x,\partial M)
> \eps
\end{equation}
(\cite{Ava:UES,Lev:UEELO} in $\R^n$, \cite{Ava:UEGRM} for manifolds,
\cite{Hor:SFEO} for higher order operators on manifolds; see 
\cite{SafVas:ADEPDO} for further improvements).
The connection to  Theorem \ref{mainthm} is established by writing
$u_{[\lam-1,\lam]}^2 = u_{[0,\lam]}^2-u_{[0,\lam-1]}^2$. Then 
\eqref{avak} gives immediately the interior estimate \eqref{int}, and
also shows the optimality of the power $(n-1)/2$ in \eqref{int}, for
{\em any} $M$ (as opposed to (\ref{int}') which is optimal only for some $M$).

It follows from the very precise and general 
results in Ivrii's book (\cite{Ivr:MAPSA}) that the $\eps$-dependence in
\eqref{avak} may be removed, and hence that our theorem even holds with $u$
replaced by $u_{[\lam-1,\lam]}$.
Our main point here is the simplification of the arguments near the boundary.
However, not all difficulties can be avoided: We still have to refer
to results on the boundary parametrix
in the treatment of the layer $\{x:\lambda^{-1}\leq \dist(x,\partial M) 
\leq 1\}$.

If one considers $\|u\|_p$ for $p\in(0,\infty)$ instead of $p=\infty$ then one
obtains interesting phenomena related to the 'restriction theorem for the
Fourier transform' of Stein-Tomas \cite{Tom:RTFT}. 
The optimal interior estimates were obtained by Sogge
\cite{Sog:CLNSCSOEOCM,Sog:FICA}. 
 The same estimate extends uniformly to
concave portions of the boundary,
 as shown in \cite{Gri:LBESPLNCB} for $n=2$
and in \cite{SmiSog:CSWEOCO} for all $n$, but
 not in general (e.g.\ for the
'whispering gallery' eigenfunctions on the
disk, see \cite{Gri:LBESPLNCB} and the remark at the end of Section
\ref{secdisk}). The problem of finding optimal $L^p$-bounds
for general boundary geometry is still open.
 
Corollary \ref{cor} also follows from the 'sharp Weyl formula' improving
\eqref{weyl}
$$
\#\{j:\,\lam_j\leq\lam\} = \gamma_M \lam^{n} + O(\lam^{n-1}).
$$
While this follows from Ivrii's results again, there are simpler proofs,
see \cite{See:ENBSFLO,Pha:MEARFSVPRL}.
\medskip

{\em Contents of the paper.} 
In Section \ref{secelem} we collect some basic facts about our problem 
which are well-known to the experts but scattered or not present in
the literature. In particular, we 
prove Corollary \ref{cor} from Theorem \ref{mainthm}
 and give two simple proofs of the interior estimate
in the case of flat domains. Also, we discuss the torus and the disk.
In Section \ref{secout} we prove  estimate \eqref{mainest}
outside a boundary layer of width $\lambda^{-1}$, and in Section \ref{secin}
we derive from this the estimate on $u$
inside this layer, for Dirichlet boundary condition.
 Finally, in Section \ref{secneumann}
we describe the modifications needed for the Neumann problem.

\section{Basic facts about $\|u\|_\infty$} \label{secelem}

\subsection{Multiplicities} 
The following proposition
 shows that Corollary \ref{cor} is  a consequence
of Theorem \ref{mainthm}.

\begin{proposition} \label{pmult}
If $V\subset L^2(M)$ is a subspace of dimension $m$,  then
\begin{equation} \label{propineq}
  \sup_{\substack{u\in V\\ \|u\|_2=1}} \|u\|_\infty
   \geq |M|^{-1/2} m^{1/2}
\end{equation}
where $|M|$ denotes the volume of $M$.
\end{proposition}

\proof
Let $v_1,\ldots,v_m$ be an orthonormal basis of $V$. For simplicity (and
sufficient for our purpose), we assume
that the $v_i$ are continuous, the general case is only slightly harder.
Define for $x,y\in M$
$$ u_y(x) = \sum_{i=1}^m \overline{v_i(y)} v_i(x),
\quad a(y) = u_y(y)= \sum_{i=1}^m |v_i(y)|^2.$$
We have $\int_M a(y) = \sum_{i=1}^m \|v_i\|_2^2 = m$, so 
$$a(\ybar)\geq m/|M|$$ for some $\ybar$.
Now by orthonormality of the $v_i$
$$ \|u_\ybar\|_2^2 = \sum_{i=1}^m |\overline{v_i(\ybar)}|^2 = a(\ybar)$$
and $\|u_\ybar\|_\infty\geq u_\ybar(\ybar) = a(\ybar)$,
so $\|u_\ybar\|^2_\infty/\|u_\ybar\|^2_2 \geq a(\ybar)\geq m/|M|$. 
Therefore, $u=u_\ybar/\|u_\ybar\|_2$
satisfies the desired bound.
\qed

\subsection{The torus} \label{sectorus}
Besides proving Corollary \ref{cor} from Theorem \ref{mainthm}, 
Proposition \ref{pmult} has another interesting consequence:
Let $M = \R^n/(2\pi\Z)^n$ be the 'square' torus. The 'standard' normalized
eigenfunctions are 
\begin{equation} \label{expfcn}
u_a(x) = (2\pi)^{-n/2} \exp(ia\cdot x),\quad a\in \Z^n,
\end{equation}
with eigenvalue $\lambda^2= |a|^2 = a_1^2 + \ldots + a_n^2$.
These are uniformly bounded, which is much better than
 the bound \eqref{thmbound}.
However, since the multiplicity of $\lambda^2$ is not uniformly bounded as
$\lambda\to\infty$, one may construct another sequence of eigenfunctions
(as in the proof of Proposition \ref{pmult}) 
with non-uniformly bounded maxima.
In fact, the multiplicity of $\lambda^2$ equals the number of representations
of $\lambda^2$ as the sum of $n$ squares of integers. From standard results
on these numbers (see \cite{Gro:RISS}, for example) one obtains:
\begin{quote}
If $M$ is the square $n$-torus, then for any $N$ 
there is a sequence of $L^2$-normalized eigenfunctions $u$ with eigenvalues $\lambda^2$ tending
to infinity and satisfying
 $$\|u\|_\infty \geq c_N \lambda^{(n-2)/2} (\log \lambda)^N$$
for some  $c_N>0$.
\end{quote}
As a simple example we take $\lambda^2 = 5^l$, $l=1,2,3,\ldots$. 
This has $4(l+1)$ representations
as sum of two squares, for example $5$ arises from $(\pm1,\pm2)$ and
$(\pm2,\pm1)$. This gives a sequence as desired for $n=2$ and $N=1$.

Note also that in the case of the torus one has equality in 
\eqref{propineq}, as follows immediately from \eqref{expfcn}.
The number-theoretic results referred to above therefore also yield
the upper bound $$ \|u\|_\infty \leq C_\epsilon \lambda^{(n-2)/2 + \epsilon}$$
for any $\epsilon>0$.
For more on the case of the torus see Bourgain's article \cite{Bou:EBLNT}.
Of course, the same results hold for a square (or cube) in $\R^n$.

See \cite{Jak:QLFT} for more about eigenfunctions on tori and
 \cite{TotZel:RMWUBE} for interesting recent work on manifolds with uniformly
bounded eigenfunctions.

\subsection{The ball} \label{secdisk}
For completeness and to show that the bound in Theorem \ref{mainthm}
cannot be improved in general, we shortly discuss the case
$M=\{x\in\R^n:\,|x|\leq 1\}$.
In polar coordinates, the Euclidean Laplacian is
\begin{equation} \label{polarcoord}
\Delta = \frac{\partial^2}{\partial r^2} + \frac{n-1}r \frac{\partial}
{\partial r} + \frac1{r^2} \Delta_S
\end{equation}
where $\Delta_S$ is the Laplacian on the unit sphere $S=\{|x|=1\}$.
By separation of variables one obtains that 
there is a basis of 
 eigenfunctions of the form
\begin{equation}\label{eqball}
u(r\omega) = r^{-(n-2)/2} J_m(\lambda r) \Phi(\omega),\quad \omega\in S,
\end{equation}
where $\Phi$ is an eigenfunction of $-\Delta_S$ (with eigenvalue $\mu^2$),
$J_m$ is the Bessel function (see \cite{Wat:TBF})
 of order $m = \sqrt{\mu^2+(n-2)^2/4}$
(which is always an integer),
and $\lambda$ is a positive zero of $J_m$ (for Dirichlet boundary conditions).

Consider a radial 
eigenfunction, i.e.\ $\mu=0$ and
 $\Phi \equiv 1$. This is not normalized as in \eqref{lzwo},
so instead of $\|u\|_\infty$ we need to estimate $\|u\|_\infty/\|u\|_2$.
We have\footnote{
That $r^{-m}J_m(r)$ attains its maximum at $r=0$ follows immediately from
Poisson's integral
$$ J_m(r) = C_m r^m \int_{-1}^1 e^{irt} (1-t^2)^{m-\frac12}\, dt,$$
see \cite{Wat:TBF}, Section 3.3. See also \cite{Wat:TBF}, Section 15.31,
for a different proof.
}
 $m=(n-2)/2$, 
$\|u\|_\infty = u(0) = C_n \lambda^{(n-2)/2}$, and the asymptotics
$$ J_m (s) \sim cs^{-1/2} \cos(s-m\pi/2-\pi/4) + O(s^{-5/2}),\quad s\to\infty$$
easily imply
$\|u\|_2^2 \approx \int_0^1 (r^{-(n-2)/2}(\lambda r)^{-1/2})^2
r^{n-1} \,dr = \lambda^{-1}$, so 
$\|u\|_\infty/\|u\|_2 \approx \lambda^{(n-1)/2}$,
and this shows that the bound \eqref{thmbound} is saturated by $u$.

Similar (but more involved)
 Bessel function estimates may be used to prove Theorem
\ref{mainthm} directly for the ball, and even the stronger estimate for
$u_{[\lam-1,\lam]}$.
See \cite{Gri:LBESPLNCB}.
\medskip

\noindent{\bf Remark:}
The radial functions are one extreme case of \eqref{eqball}.  The other
extreme case is obtained by taking $\lambda=\lambda_{m1}$, the first positive
zero of $J_m$. For $m=1,2,3,\ldots$ this yields the sequence of 'whispering
gallery eigenfunctions' (say $n=2$ for simplicity). 
They concentrate on a strip of width $\approx
\lambda^{-2/3}$ at the boundary, as follows from (and is made precise by)
the estimates
\begin{align*}
\lambda_{m1} &= m + am^{1/3} + O(m^{-2/3}),\quad \text{as }m\to\infty,
                \text { with } a>0\\
J_m (m+tm^{1/3}) &\geq m^{-1/3} \quad\text{ for } t\in[-a/2,a/2]\\
J_m(m+tm^{1/3}) &\leq Cm^{-1/3} e^{-c|t|^{3/2}}\quad\text{ for }
-m^{2/3}\leq t\leq 2a
\end{align*}
for positive constants $c,C$.
(These are easy consequences of well-known asymptotic formulas
for $J_m$, see 
\cite{Olv:AEBFLO} for example, or \cite{Wat:TBF}, Sec.8.4 for weaker but 
sufficient bounds.)
Sogge showed that away from the boundary  concentration can happen only
on  sets of area $\geq \lambda^{-1/2}$.
This follows from his estimate $\|u\|_6 \leq C\lambda^{1/6}\|u\|_2$.
In contrast, the whispering gallery eigenfunctions have 
$\|u\|_6\approx \lambda^{1/3}\|u\|_2$.
Note that the $L^\infty$ estimate \eqref{mainest} only implies area $\geq
\lambda^{-1}$ for concentration. This shows two things:
\begin{enumerate}
\item
Optimal bounds on concentration phenomena are obtained from certain $L^p$,
$p<\infty$,
rather than $L^\infty$ bounds on eigenfunctions.
\item
As opposed to $L^\infty$ bounds, these $L^p$-bounds are sensitive to the 
presence (and geometry) of a boundary.
\end{enumerate}
In general, such optimal bounds are still unknown.

\subsection{General domains in $\R^n$} \label{seceucl}
Here we prove the interior estimate \eqref{int} for Euclidean
domains $M\subset\R^n$. 
We give two proofs: First, we finish the wave equation
proof outlined in the Introduction, and second, we give a more
direct proof using averaging and Bessel functions.

{\em First proof:}
\newcommand{\KR}{K_{\R^n}}
As argued in the introduction, it is sufficient to prove \eqref{mainest},
with $K$ replaced by the wave kernel in $\R^n$. One way to represent
$\KR$ is as an oscillatory integral, obtained from using the
$x$-space Fourier transform
and solving the  ordinary differential equation  that results from
\eqref{waveeqn},
\begin{equation} \label{Kformula}
\KR(t,x,y) = (2\pi)^{-n} \int_{\R^n} \cos (t|\xi|)\, e^{i(x-y)\xi}\,d\xi.
\end{equation}
From this it is not hard to get \eqref{mainest} directly. Let us
describe a more conceptual approach which yields the power $n-1$ by pure
homogeneity arguments (this is the scaling technique in
\cite{Mel:TWG}
 for this special case): 
From the homogeneity of the equation satisfied by $\KR$ (or from
\eqref{Kformula} directly) one has the homogeneity
$$ \KR(t,x,y) = \eps^{-n} \KR(t/\eps,x/\eps,y/\eps)$$
for $\eps>0$. Since $\Delta$ has constant coefficients, one has furthermore
translation invariance $\KR(t,x,y) = \KR(t,x-y,0)$. Therefore,
$\KR(t,x,x)=\KR(t,0,0)$ is a distribution in $t$, homogeneous of degree $-n$.
Then its  singular support must be contained in  $\{0\}$, and
its inverse Fourier transform is also smooth
outside zero and homogeneous of degree $n-1$ (see \cite{Hor:ALPDOI}, Theorem 7.1.18), 
so it 
satisfies the desired bound. Then it is straightforward to see that
$\rho * \KR\spcheck$ satisfies the same bound, i.e. \eqref{mainest}.

{\em Second proof:} 
Let $x_0\in M$, and assume that the ball $B$ of radius
$R$ around $x_0$ is contained in $M$. We show that
\begin{equation} \label{localest}
|u(x_0)| \leq C\lambda^{(n-1)/2} R^{-1/2} \|u\|_{L^2(B)}
\end{equation}
which clearly implies (\ref{int}'). 
To simplify notation, assume $x_0=0$. In this proof, $\|\,\|_p$ denotes
the $L^p$ norm on $B$, not $M$.

Define the spherical average 
$$h(r) = \frac1{|S|} \int_S u(r\omega)\,d\omega, $$
considered as function on $B$.
Since $h$ is the  average over the functions $g^*u$
over all rotations $g\in SO(n)$ and since these rotations induce
isometries on $L^2(B)$, Minkowski's inequality gives 
$$\|h\|_2 \leq \|u\|_2.$$
Furthermore, $h$ solves $(\Delta+\lambda^2)h=0$ (since $\Delta$ is 
rotation invariant) and is radial, so it
is of the form \eqref{eqball} with $\Phi=const$. Therefore,
the same calculation as in Section \ref{secdisk} shows that
$$\frac{\|h\|_\infty}{\|h\|_2} \approx R^{-1/2}\lambda^{(n-1)/2}.$$
Finally, $\|h\|_\infty = |h(0)|=|u(0)|$, so \eqref{localest} follows.
We remark that this proof may be adapted to prove \eqref{int} instead
of (\ref{int}'), see \cite{Gri:LBESPLNCB}.

\section{Estimates outside the boundary layer}  \label{secout}
In this section we prove Theorem \ref{mainthm}\ for points $x$ with
$$ \dist(x,\partial M) \geq \lam^{-1}.$$
As explained in the Introduction, we only have to prove (\ref{mainest})
for these $x$.
Propagation of
singularities (see \cite{Hor:ALPDOIII})
implies that the singular support of the distribution $t\ra K(t,x,x)$ is
contained in the set of lengths of geodesics, possibly 
reflected at the boundary,  
which start and end at $x$. Clearly, for small $\epsilon$ the only 
singularities in $|t|<\eps$ are therefore at $t=0$ 
and possibly at $t=\pm 2\dist(x,\partial M)$.
Therefore, $K$ may be expected to be  and indeed is 
representable, for small $|t|$, as the
sum of two distributions, a 'direct' term which is only singular at $t=0$,
and a 'reflected' term. To describe the direct term,
choose a closed manifold $\mtilde$ extending $M$, and let 
$$\kdir(t,x,y)$$ be the solution of the problem (\ref{waveeqn}) on $\mtilde$.

To describe the reflected term, it is convenient to introduce geodesic
coordinates with respect to the boundary; that is, we identify points $x$
of $M$ close to the boundary with pairs 
$(x',x_n)\in\partial M\times [0,c)$, $c>0$, via the map
\begin{equation} \label{normalcoord}
\partial M \times [0,c) \ra M
\end{equation}
sending $(x',x_n)$ to the endpoint of the geodesic of length $x_n$ which
starts at $x'\in\partial M$ perpendicular to $\partial M$. 
This map is
a diffeomorphism onto its image for sufficiently small $c$.
We have
$$ x_n = \dist (x,\partial M).$$
Also, let $\chi_+^\alpha$, $\alpha\in\C$, be the distribution on $\real$ obtained by 
analytic continuation in $\alpha$ from $\{\text{Re } \alpha > -1\}$, where
it is defined by
\begin{equation*}
\chi_+^\alpha (s) = \begin{cases} s^\alpha/\Gamma(\alpha+1) & \text{ if }
				       s > 0\\
				    0 & \text{ if } s\leq 0
                    \end{cases}
\end{equation*}
and satisfies $(\chi_+^{\alpha+1})' = \chi_+^\alpha$, see \cite{Hor:ALPDOI}, Section 3.2.
Alternatively, it is the inverse Fourier transform of 
\begin{equation} \label{chift}
\hat{\chi}_+^\alpha (\sigma) = e^{-i(\alpha+1)\pi/2}
			\frac1{(\sigma - i0)^{\alpha+1}}.
\end{equation}

Theorem 17.5.9 in \cite{Hor:ALPDOIII} then gives the following description of the
singularity of $K$ at $t=\pm 2 \dist(x,\partial M)$:

\begin{proposition} \label{parametrix}
For sufficiently small $\eps$, there is a distribution
$$I(x',\theta,t) \in \distr (\partial M \times \real
				  \times (-\eps,\eps)) $$
so that for $|t|<\eps$ and $x_n>0$ we have
$$K(t,x,x) = \kdir(t,x,x) - t^{-n} I(x',\frac{2x_n}t,t).$$
Furthermore, $I$ has support in $|\theta|\leq1$ and singular support in $|\theta|=1$,
and
near $\theta = 1$ we have 
$$I(x',\theta,t) = \sum_{j=0}^{N-1} a_j(x',\theta,t) 
		      \chi_+^{j-(n+1)/2}(1-\theta) + R_N(x',\theta,t)$$
with smooth functions $a_j$ and a continuous remainder $R_N$, for any
$N>(n-1)/2$. A similar expansion exists near $\theta=-1$.
\end{proposition}

Note: In the statement of Theorem 17.5.9 in \cite{Hor:ALPDOIII} the direct term
(called $t^{-n}I_1$ there) is not identified like in the statement above.
But in the proof (middle of page 59, loc.\ cit.) 
it is chosen as $\kdir$ like above.

The proposition is obtained by analyzing the regions $t<3x_n$
and $t\geq 3x_n$ separately: In $t<3x_n$ one uses a 
suitable scaling (cf. Section \ref{seceucl}) and the Hadamard parametrix,
and in $t\geq3x_n$ (where no singularities should occur) one uses 
propagation of singularities estimates (this is the harder part).

Since the estimates for $\kdir$ are already known by the interior estimate, 
(\ref{mainest}) will
follow from:

\begin{lemma}
For $x_n\geq \lam^{-1}$, we have
$$ |\int_0^\infty e^{it\lam} \rhohat(t) t^{-n}I(x',\frac{2x_n}t,t)\, dt|
	  \leq C\lam^{n-1}.$$
\end{lemma}
For convenience, we restrict to positive $t$. Negative $t$ are handled
in the same way.

\proof
The integrand is zero for $t<2x_n$ and for $t>\epsilon$. 
We split the integral up
into a part where $t>3x_n$ and a part where $t<4x_n$, using a cutoff
function smooth in $t/x_n$.
By Proposition \ref{parametrix}, 
$I$ is smooth on the first part and therefore bounded,
so this part is dominated by a constant times
$$ \int_{x_n}^\infty t^{-n}\,dt \leq x_n^{-n+1} \leq \lam^{n-1}.$$
In the second part, where $t<4x_n$, 
we split up $I$ as in the proposition, with a fixed
$N>(n-1)/2$. Then the
term $R_N$ can be handled in the same way as the first part. 

It remains to analyze the singular terms. 
If we denote the cutoff function by $\psi(t/x_n)$, $\psi\in C_0^\infty (1,4)$, they
take the form (assuming for simplicity that $\rhohat$ was chosen
to be constant near zero)
$$\int e^{it\lam} a(x',\frac{2x_n}t,t) \chi_+^\alpha(1-\frac{2x_n}t)
\psi(\frac{t}{x_n})t^{-n}\, dt, $$ 
$\alpha = j-(n+1)/2$.
If we change variables $\tau=t/2x_n$ and use the homogeneity of $\chi_+^\alpha$,
this becomes
$$x_n^{-n+1} \int e^{2i\tau \lam x_n} b(x',\tau,t) 
		      \chi_+^\alpha (\tau-1)\,d\tau$$
with $b$ smooth and supported in $\tau\in(1/2,2)$.
Using (\ref{chift}) one sees
by a short standard calculation that this is bounded by a constant times
$$ x_n^{-n+1} (\lam x_n)^{(n-1)/2 - j} 
      = \lam^{n-1} (\lam x_n)^{-(n-1)/2-j} \leq \lam^{n-1}.$$
\qed

\section{Estimates in the boundary layer} \label{secin}
So far, we have proved $|u(x)| \leq C\lam^{(n-1)/2}$ for $\dist(x,\partial M)
 \geq\lam^{-1}$.
The proof of Theorem \ref{mainthm}\ in the Dirichlet case
 will therefore be completed 
by the following lemma:

\begin{lemma} \label{lmaxpr}
If $u$ is a solution of $$(\Delta+\lam^2)u=0$$ vanishing at the
boundary of $M$, then 
\begin{equation}\label{eqmaxpr}
\max_{x:\dist(x,\partial M)<\lamone} |u(x)| \leq
	    \max_{x:\dist(x,\partial M)= \lamone} |u(x)|. 
\end{equation}
\end{lemma}

\proof
Without loss of generality, assume $u$ is real valued.
Under the identification near the boundary
given by (\ref{normalcoord}), the metric takes the
form
$$ g(x) = dx_n^2 + g'(x',x_n)$$
where $g'(\cdot,x_n)$ is a Riemannian metric on $\partial M$ for each $x_n$.
Therefore, the Laplacian has the form
\begin{equation} \label{laplacian}
\Delta = \frac{\partial^2}{\partial x_n^2} + 
      a(x)\frac{\partial}{\partial x_n}+ P(x',x_n,D_{x'})
\end{equation}
with a Laplacian $P$ on the boundary (depending on the parameter
$x_n$) and a smooth function $a$.
The generalized maximum principle
says that if $v$ is a positive
function on the strip
$$ S = \{0\leq x_n \leq \lamone\} \subset M $$
with $$(\Delta  + \lam^2)v \leq 0 $$
then 
\begin{equation}\label{maxpr}
 \max_S \frac{|u|}v \leq \max_{\partial S} \frac{|u|}v.
\end{equation}
(Apply \cite{ProWei:MPDE}, Theorem 10, to $u$ and then to $-u$.)
We apply this with the function
$$ v(x',x_n) = \sin (\frac\pi2 + \frac32(\lam x_n-1)).$$
We have $v>0.07$ on $S$ and thus, from (\ref{laplacian}),
$$ (\Delta +\lam^2)v = -\frac54 \lam^2 v + O(\lam) < 0 \text{ on $S$
 for large } \lam,$$
so the generalized maximum principle applies.
Now $u=0$ at the outer boundary $x_n=0$ of $S$, and $v=1$ at the inner
boundary $x_n=\lam^{-1}$. Since $v\leq 1$ on $S$, \eqref{maxpr} implies
\eqref{eqmaxpr}.
\qed

\section{The Neumann problem} \label{secneumann}
Theorem \ref{mainthm} and Corollary \ref{cor} remain true if the Dirichlet
boundary condition $u_{|\partial M}=0$ is replaced by the Neumann
boundary condition $$\partial_n u_{|\partial M} = 0$$ in \eqref{evproblem},
where $\partial_n$ denotes the outward normal derivative.

Let us sketch the proof: 
The proof of Proposition \ref{pmult} (and thus the proof of Corollary 
\ref{cor} from Theorem \ref{mainthm}) as well as the reduction of 
Theorem \ref{mainthm} to the wave kernel estimate \eqref{mainest}
carry over literally, except that $K$ is replaced by the Neumann
wave kernel (i.e. $\partial_nK_{|x\in\partial M}=0$ instead of
$K_{|\partial M}=0$ in \eqref{waveeqn}, where $\partial_n$ refers to
the $x$-coordinates).
(Also, Section \ref{seceucl} applies literally to the Neumann problem.)

Proposition \ref{parametrix} holds for the Neumann wave kernel as well
(with different $a_j$), by straightforward modification of the arguments
in \cite{Hor:ALPDOIII}. 

Finally, Lemma \ref{lmaxpr} remains true for Neumann eigenfunctions, except
that \eqref{eqmaxpr} must be replaced by
\begin{equation} \tag{\ref{eqmaxpr}'}
\max_{x:\dist(x,\partial M)<\lamone} |u(x)| \leq
	    20\,\max_{x:\dist(x,\partial M)= \lamone} |u(x)|. 
\end{equation}
To see this, just note that, in the comparison of the functions
$u$ and $v$, the maximum of $u/v$ over $S$ not only must occur at a boundary
point $x_0$ of $S$, 
but also the outward normal derivative at $x_0$ must be strictly positive.
See Theorem 10 of \cite{ProWei:MPDE} (this is sometimes called Zaremba's
principle).
Now choose $v = \sin(\frac\pi2 + \frac32 \lam x_n)$, then
$$ \partial_n v = -\partial v/\partial x_n = 0\quad\text{ for } x_n = 0$$
and so
$\partial_n (u/v) = 0$ for $x_n=0$.
 This means that the maximum
of $u/v$ over  $S$ must occur at the inner boundary
$x_n=\lambda^{-1}$ of $S$.
Since $v^{-1}<20$ there, we get (\ref{eqmaxpr}') 
after applying the same argument to $-u$.


\end{document}